\documentclass[12pt,leqno]{amsart}
\usepackage{amsmath}
\usepackage{amssymb}

\def\underset#1#2{{\mathrel{\mathop {{}_{} {#2}}\limits_{{#1}_{}}}}}
\def\upplim_#1{\underset{#1}{\overline\lim}\;}
\def\lowlim_#1{\underset{#1}{\underline\lim}\;}

\def\le{\leqslant}
\newtheorem{corollary}[equation]{Corollary}
\newtheorem{definition}[equation]{\indent{\it Definition}\rm }

\newtheorem{lem}[equation]{Lemma}
\newtheorem{proposition}[equation]{Proposition}

\newtheorem{theorem}[equation]{Theorem}

\newcommand{\C}{{\mathbf{C}}}

\renewcommand{\P}{{\mathbf{P}}}

\newcommand{\pnc}{{\mathbf{P}^n{(\mathbf{C})}}}

\numberwithin{equation}{section}

\title[{\tiny Uniqueness problem of meromorphic mappings from a complete K\"{a}hler manifold}]{Uniqueness problem of meromorphic mappings from a complete K\"{a}hler manifold into a projective variety}

\begin{document}

\date{} 
\author{Le Ngoc Quynh}
\address{Faculty of Education\\
An Giang University\\
18 Ung Van Khiem, Dong Xuyen, Long Xuyen, An Giang, Vietnam}
\email{nquynh1511@gmail.com}

\thanks{This research is funded by Vietnam National Foundation for Science and Technology Development (NAFOSTED) under grant number 101.04-2015.03.}


\maketitle

\begin{abstract} The purpose of this paper has twofold. The first is to prove a unicity theorem for meromorphic mappings of a complete K\"{a}hler manifold $M$ in $\pnc$ sharing few hypersurfaces. The second is to give a unicity theorem for the case of differential nondegenerate meromorphic mappings sharing $2N-n+3+2\rho(2N-n+1)$ hyperplanes in $N-$subgeneral position.
\end{abstract}

\section{Introduction}

In 1926, Nevanlinna \cite{N} showed that, two distinct non-constant meromorphic functions $f$ and $g$ on the complex plane $\mathbb{C}$ can not have the same inverse images for five distinct values. Over the last few decades, there have been many generalizations of Nevanlinna's result to the case of meromorphic maps from $\mathbb{C}^m$ into the complex projective space $\mathbb{P}^n(\mathbb{C})$.

In \cite{Fu2}, Fujimoto gave the following new type of uniqueness theorem for meromorphic maps of an $m$-dimensional complete K\"{a}hler manifold $M$ into $\mathbb{P}^n(\mathbb{C})$, where the universal covering of $M$ is biholomorphic to the ball in $\mathbb{C}^m.$

\vskip0.2cm
\noindent
\textbf{Theorem A.} {\it
Let $M$ be a complete, connected K\"{a}hler manifold whose universal is biholomorphic to $B(R_0)\subset \mathbb{C}^m$, where $0<R_0\le \infty.$ Let $f$ and $g$ be linearly non-degenerate meromorphic maps of $M$ into $\mathbb{P}^n(\mathbb{C})$ and $H_1,\ldots,H_q$ be $q$ hyperplanes in general position with $\dim f^{-1}(H_i\cap H_j)\le m-2$ for $1\le i<j\le q.$ Assume that $f$ and $g$ satisfy the condition $(C_\rho)$ for some non-negative constant $\rho$, and

(i) $f^{-1}(H_i)=g^{-1}(H_j)$ for $1\le j\le q,$

(ii) $f=g$ on $\cup_{j=1}^q f^{-1}(H_j).$

\noindent
If $q>n+1+\rho n(n+1),$ then $f\equiv g.$}

\vskip0.2cm
Here, we say that the mapping $f$ satisfies the condition $(C_\rho)$ for some non-negative constant $\rho$ if there exists a non-zero bounded continuous real-valued function $h$ on $M$ such that
$$\rho \Omega_f + dd^c\log h^2\ge \text{Ric}\omega,$$
where $\Omega$ is the Fubini-Study form on $\mathbb{P}^n(\mathbb{C})$, $\omega=\frac{\sqrt{-1}}{2}\sum_{i,j}h_{i\overline{j}}dz_i\wedge dz_{\overline{j}}$ is K\"{a}hler form on $M$, $\text{Ric}\omega=dd^c\log(det(h_{i\overline{j}}))$, $d=\partial+\overline{\partial}$ and $d^c=\frac{\sqrt{-1}}{4\pi}(\overline{\partial}-\partial)$

Recently, Ru and Sogome \cite{R-S2} have extended this result to the case where the meromorphic maps share hypersurfaces instead of hyperplanes. In this paper, we will improve and extend the result of Ru-Sogome by giving  a uniqueness theorem with less than involving hypersurfaces and the hypersurfaces may be located in $N$-subgeneral positon.

To state our results, we recall the following notation due to Thai - Quang \cite{T-Q}. 

Let $N\ge n$ and $q\ge N+1$. Let $Q_1,\ldots,Q_q$ be hypersurfaces in $\mathbb{P}^n(\mathbb{C}).$ The hypersurfaces $Q_1,\ldots,Q_q$ are said to be in $N$-subgeneral position in $\mathbb{P}^n(\mathbb{C})$ if $Q_{i_0}\cap\cdots\cap Q_{i_N}=\emptyset$ for every $1\le i_0<\cdots<i_N\le q.$

Let $V$ be a complex projective subvariety of $\mathbb{P}^n(\mathbb{C})$ of dimension $k$ $(k\le n)$. Let $Q_1,\ldots,Q_q$ $(q\ge k+1)$ be $q$ hupersurfaces in $\mathbb{P}^n(\mathbb{C})$. The family of hypersurfaces $\{Q_i\}_{i=1}^q$ is said to be in $N$-subgeneral position with respect to $V$ if for any $1\le i_1<\cdots<i_{N+1}\le q,$ 
$$V\cap (\bigcup\limits_{j=1}^{N+1}Q_{i_j})=\emptyset.$$
If $\{Q_i\}_{i=1}^q$ is in $n$-subgeneral position then we say that it is in general position with respect to $V$.

Now, let $d$ be a positive integer. We denote by $I(V)$ the ideal of homogeneous polynomials in $\mathbb{C}[x_0,\ldots,x_n]$ defining $V$ and by $H_d$ the vector space of all homogeneous polynomials in $\mathbb{C}[x_0,\ldots,x_n]$ of degree $d$. Define
$$I_d(V):=\frac{H_d}{I(V)\cap H_d} \ \text{and} \ H_V(d):=\dim I_d(V).$$
Then $H_V(d)$ is called the Hilbert function of $V$. Each elements of $I_d(V)$ which is an equivalent class of an elements $Q\in H_d$, will be denoted by $[Q].$

Let $M$ be a complete K\"{a}hler manifold. Let $f:\mathbb{C}^m\to V$ be a meromorphic mapping. We say that $f$ is degenerate over $I_d(V)$ if there is $[Q]\in I_d(V)\backslash\{0\}$ such that $Q(f)=0$. Otherwise, we say that $f$ is nondegenerate over $I_d(V)$. It is clear that if $f$ is algebraically nondegenerate then $f$ is nondegenerate over $I_d(V)$ for every $d\ge 1$. 

In \cite{T-Q}, Thai - Quang gave a non-integrated defect relation as follows.

\vskip0.2cm
\noindent
\textbf{Theorem B.} {\it 
Let $V$ be a complex projective subvariety of $\mathbb{P}^n(\mathbb{C})$ of dimension $k$ $(k\le n)$. Let $\{Q_i\}_{i=1}^q$ be hypersurfaces of $\mathbb{P}^n(\mathbb{C})$ in N-subgeneral position with respect to $V$ with $degQ_i=d_i \ (1\le i\le q)$. Let $d$ be a least common multiple of $d'_is$, i.e., $d=lcm(d_1,\ldots,d_q).$ Let $f$ be a meromorphic mapping of a complete K\"{a}hler manifold $M$, whose universal covering is biholomorphic to $\mathbb{C}^m$ or a ball $B(R) \subset \mathbb{C}^m$, into $V$ such that $f$ is nondegenerate over $I_d(V)$. Assume that there is a continuous plurisubharmonic function $h$ on $M$, a position constant $\rho$ such that
$$\rho f^*\Omega + \frac{\sqrt{-1}}{2\pi}\partial \overline{\partial}\log h^2 \le Ric \omega,$$
where $\omega$ is Kaler form on $M$ and $\Omega$ is the Fubini-Study form on $\mathbb{P}^n(\mathbb{C})$. Then, we have
$$\sum_{j=1}^q\delta_{H_V(d)-1}^f(Q_j) \le \frac{(2N-k+1)H_V(d)}{k+1}+\frac{\rho H_V(d)(H_V(d)-1)}{d}.$$	}

Actually, we may assume that $M=B(R_0)(\subset\mathbb{C}^m)$, where $0<R_0\le\infty.$ Using some idea from the proof of Thai and Quang, in this paper, we will prove a uniqueness theorem as follows.

\begin{theorem}\label{1.1}
 Let $V$ be a complex projective subvariety of $\mathbb{P}^n(\mathbb{C})$ of dimension $k \ (k\le n).$ Let $\{Q_i\}_{i=1}^{q}$ be hypersurfaces of $\mathbb{P}^n(\mathbb{C})$ in $N$-subgeneral position with respect to $V$ with degree of $d_i \ (1\le i\le q)$. Let $d$ be the least common multiple of $d'_i s,i.e,d=lcm(d_1,\ldots,d_q).$ Let $f, g: B(R_0) (\subset\mathbb{C}^m) \to V$ be nondegenerate over $I_d(V)$ such that $\lim \limits_{r\to R_0}\sup\frac{T_f(r,r_0)+T_g(r,r_0)}{\log (1/(R_0-r))}<\infty$, where $0<r_0<r<R_0<\infty.$ Assume that $f$ and $g$ satisfy the following conditions

(i)\ $f=g$ on $\bigcup\limits_{i=1}\limits^{q}(Zero Q_i(f)\cup Zero Q_i(g)),$

(ii)\ there exists a non-negative constant $\rho$ such that $f$ and $g$ satisfy the condition $(C_\rho).$

\noindent
If $q>\frac{(2N-k+1)H_V(d)}{k+1}\bigl(1+\frac{\rho(H_V(d)-1)}{d}\bigl)+\frac{2N(H_V(d)-1)}{d},$ then $f\equiv g.$ 
\end{theorem}

For the case where $R_0=\infty$ or $\lim \limits_{r\to R_0}\sup\frac{T_f(r,r_0)+T_g(r,r_0)}{\log (1/(R_0-r))}=\infty$, we have the following corollary.

\begin{corollary}\label{1.2}
 Let $V$ be a complex projective subvariety of $\mathbb{P}^n(\mathbb{C})$ of dimension $k \ (k\le n).$ Let $\{D_i\}_{i=1}^{q}$ be hypersurfaces of $\mathbb{P}^n(\mathbb{C})$ in $N$-subgeneral position with respect to $V$ with degree of $d_i \ (1\le i\le q)$. Let $d$ be the least common multiple of $d'_i s,$ i.e., $d=lcm(d_1,\ldots,d_q).$ Let $f, g: B(R_0) (\subset\mathbb{C}^m) \to V$ be nondegenerate over $I_d(V).$ Assume that
 
 (i)\ $f=g$ on $\bigcup\limits_{i=1}\limits^{q}(Zero Q_i(f)\cup Zero Q_i(g)),$
 
 (ii)\ $q>\frac{(2N-k+1)H_V(d)}{k+1}+\frac{2N(H_V(d)-1)}{d},$

\noindent 
If $R_0=\infty$ or $\lim \sup \limits_{r\to R_0}\frac{T_f(r,r_0)+T_g(r,r_0)}{\log (1/(R_0-r))}=\infty$,  then $f\equiv g.$ 
\end{corollary}

On the other hand, in 1981, Drouilhet \cite{D} proved a uniqueness theorem for differential nondegenerate meromorphic mappings as follows.

\vskip0.2cm 
\noindent
\textbf{Theorem C.} {\it Let $M$ be an $m$ - dimensional smooth affine algebraic variety. Let $f,g: M\to\mathbb{P}^n(\mathbb{C}) $ be differential nondegenerate meromorphic maps, $mge n.$ Let $A$ be a hypersurface of degree at least $n+4$ in $\mathbb{P}^n(\mathbb{C})$ having normal crossing. Suppose $f^{-1}(A)=g^{-1}(A)$ as point sets and $f$ and $g$ agree at all points of $f^{-1}(A)$ lying in their common domain of determinacy. Suppose either $M=\mathbb{C}^m$ or $f$ and $g$ are transcendental. Then $f=g$.} 

By using technique in \cite{T-Q} , we will generalize the above result of Drouilhet for the case where $M$ is a K\"{a}hler manifold and the targets are hyperplanes in subgeneral position. Firstly, we recall the defect relation given by Thai - Quang in \cite{T-Q} as follows.

\vskip0.2cm 
\noindent
\textbf{Theorem D.} {\it Let $\{H_i\}_{i=1}^q$ be hyperplanes of $\mathbb{P}^n(\mathbb{C})$ in N-subgeneral position. Let $f$ be a meromorphic mapping of complete Kahler manifold $M$ of dimension $m$, $(m \ge n)$, whose universal covering is biholomorphic to $\mathbb{C}^m$ or a ball $B(R) \subset \mathbb{C}^m$, into $\mathbb{P}^n(\mathbb{C})$ such that $f$ is differential nondegenerate. Assume that there is a continuous plurisubharmonic function $h$ on $M$, a position constant $\rho$ such that
	$$\rho f^*\Omega + \frac{\sqrt{-1}}{2\pi}\partial \overline{\partial}\log h^2 \le Ric \omega,$$
	where $\omega$ is Kaler form on $M$ and $\Omega$ is the Fubini-Study form on $\mathbb{P}^n(\mathbb{C})$. Then, we have
	$$\sum_{j=1}^q\delta_1^f(H_j) \le 2N-n+1+2\rho(2N-n+1).$$}
Hence, in the above theorem, if the family of hyperplanes is assumed to be in general position, the defect relation will be obtained as follows
$$\sum_{j=1}^q\delta_1^f(H_j) \le n+1+2\rho(n+1)$$

By the condition that the universal covering of $M$ is biholomorphic to $\mathbb{C}^m$ or a ball $B(R) \subset \mathbb{C}^m$, we may assume that $M=B(R_0)(\subset\mathbb{C}^m)$, where $0<R_0\le\infty$, and prove the following.

\begin{theorem}\label{1.3}
	Let $\{H_i\}_{i=1}^{q}$ be hyperplanes of  $\mathbb{P}^n(\mathbb{C})$ in $N$-subgeneral position. Let $f, g: B(R_0)(\subset\mathbb{C}^m) \to \mathbb{P}^n(\mathbb{C}) \ (n\le m)$ be meromorphic mappings such that $f$ and $g$ are differential nondegenerate and $\lim \limits_{r\to R_0}\sup\frac{T_f(r,r_0)+T_g(r,r_0)}{\log (1/(R_0-r))}<\infty$, where $0<r_0<r<R_0<\infty$. Assume that $f$ and $g$ satisfy the following conditions
	
	(i)\ $f^{-1}(H_i)=g^{-1}(H_i) \ \text{for} \ 1\le i\le q,$ 
	
	(ii)\ $f=g$ on $\bigcup_{i=1}^q f^{-1}(H_i),$
	
	(iii)\ there exists a non-negative constant $\rho$ such that $f$ and $g$ satisfy the conditon $(C_\rho).$
	
\noindent
	If $q>2N-n+3+2\rho(2N-n+1),$ then $f\equiv g.$ 
\end{theorem}

For the case where $R_0=\infty$ or $\lim \limits_{r\to R_0}\sup \frac{T_f(r,r_0)+T_g(r,r_0)}{\log (1/(R_0-r))}=\infty$, we have the following corollary.

\begin{corollary}\label{1.4}
	Let $\{H_i\}_{i=1}^{q}$ be hyperplanes of $\mathbb{P}^n(\mathbb{C})$ in $N$-subgeneral position. Let $f, g: M \to \mathbb{P}^n(\mathbb{C}) \ (n\le m)$ be differential nondegenerate. Assume that $f$ and $g$ satisfy the following conditions
	
		(i)\ $f^{-1}(H_i)=g^{-1}(H_i) \ \text{for} 1\le i\le q,$ 
		
		(ii)\ $f=g$ on $\bigcup_{i=1}^q f^{-1}(H_i),$
		
		(iii)\ $q>2N-n+3.$ 
	
\noindent
	If $R_0=\infty$ or $\lim \limits_{r\to R_0}\sup \frac{T_f(r,r_0)+T_g(r,r_0)}{\log (1/(R_0-r))}=\infty$,  then $f\equiv g.$ 
\end{corollary}

Hence, when $\{H_i\}_{i=1}^q$ locate in general position, we get again the result of Drouilhet \cite{D}.

\section{Basic notions and auxiliary results from Nevanlinna theory}

\noindent
{\bf (a)}\ We set $||z|| = \big(|z_1|^2 + \dots + |z_n|^2\big)^{1/2}$ for
$z = (z_1,\dots,z_n) \in \mathbb{C}^m$ and define
\begin{align*}
B(r) := \{ z \in \mathbb{C}^m : ||z|| < r\},\quad
S(r) := \{ z \in \mathbb{C}^m : ||z|| = r\}\ (0<r<\infty).
\end{align*}

Define 
$$v_{m-1}(z) := \big(dd^c ||z||^2\big)^{m-1}\quad \quad \text{and}$$
$$\sigma_m(z):= d^c \text{log}||z||^2 \land \big(dd^c \text{log}||z||^2\big)^{m-1}
\text{on} \quad \mathbb{C}^m \setminus \{0\}.$$

\noindent{\bf (b)}\ Let $F$ be a nonzero holomorphic function on a domain $\Omega$ in $\mathbb{C}^m$. For a set $\alpha = (\alpha_1,...,\alpha_m) $ of nonnegative integers, we set $|\alpha|=\alpha_1+...+\alpha_m$ and 
$\mathcal {D}^\alpha F=\frac {\partial ^{|\alpha|} F}{\partial ^{\alpha_1}z_1...\partial ^{\alpha_m}z_m}.$
We define the map $\nu_F : \Omega \to \mathbb{Z}$ by
$$\nu_F(z):=\max\ \{l: \mathcal {D}^\alpha F(z)=0 \text { for all } \alpha \text { with }|\alpha|<l\}\ (z\in \Omega).$$

We mean by a divisor on a domain $\Omega$ in $\mathbb{C}^m$ a map $\nu : \Omega \to \mathbb{Z}$ such that, for each $a\in \Omega$, there are nonzero holomorphic functions $F$ and $G$ on a connected neighbourhood $U\subset \Omega$ of $a$ such that $\nu (z)= \nu_F(z)-\nu_G(z)$ for each $z\in U$ outside an analytic set of dimension $\le m-2$. Two divisors are regarded as the same if they are identical outside an analytic set of dimension $\le m-2$. For a divisor $\nu$ on $\Omega$ we set $|\nu| := \overline {\{z:\nu(z)\ne 0\}},$ which
is a purely $(n-1)$-dimensional analytic subset of $\Omega$ or empty.

Take a nonzero meromorphic function $\varphi$ on a domain $\Omega$ in $\mathbb{C}^m$. For each $a\in \Omega$, we choose nonzero holomorphic functions $F$ and $G$ on a neighbourhood $U\subset \Omega$ such that $\varphi = \frac {F}{G}$ on $U$ and $\dim (F^{-1}(0)\cap G^{-1}(0))\le m-2,$ and  we define the divisors $\nu^0_\varphi,\ \nu^\infty_\varphi$ by 
$ \nu_\varphi := \nu_F,\  \nu^\infty_\varphi :=\nu_G$, 
which are independent of choices of $F$ and $G$ and so globally well-defined on $\Omega$. 

\noindent{\bf (c)} For a divisor $\nu$ on $\mathbb{C}^m$ and for $M$ be a positive integer or $+\infty$. We define the truncated divisor $\nu^{[M]}$ by
	$$\nu^{[M]}(z):= \min\{\nu(z), M\}$$
and the counting function of $\nu$ by
\begin{align*}
	N(r,r_0,\nu) := \int\limits_{r_0}^r \frac{n(t)}{t^{2m-1}}dt\quad
	(1 < r_0<r< R),
\end{align*}
where
\begin{align*}
	n(t): =
	\begin{cases}
		\int\limits_{|\nu| \cap B(t)} \nu(z) v_{m-1} &\text{ if } m \geq 2,\\
		\sum_{|z| \le t} \nu(z) & \text{ if } m = 1.
	\end{cases}
\end{align*}
Similarly, we define $n^{[M]}(t)$ and define $N(r,r_0,\nu^{[M]})$ which denotes by $N^{[M]}(r,r_0,\nu)$.

We omit the character $^{[M]}$ if $M=+\infty$.

\noindent{\bf (d)}\ Let $f : B(R)\subset\mathbb{C}^m \longrightarrow \mathbb{P}^n(\mathbb{C})$ be a meromorphic mapping. We say that $f$ is differential  nondegenerate if $df$ has the maximal rank. For arbitrarily fixed homogeneous coordinates $(w_0 : \dots : w_n)$ on $\mathbb{P}^n(\mathbb{C})$, we take a reduced representation $ (f_0 : \dots : f_n)$ of $f$. This means $f(z) = \big(f_0(z) : \dots : f_n(z)\big)$ outside the analytic set $\{ f_0 = \dots = f_n= 0\}$ of codimension $\geq 2$. Set $\Vert f \Vert = \big(|f_0|^2 + \dots + |f_n|^2\big)^{1/2}$.

The characteristic function of $f$ (with respect to $\Omega$) is defined by
\begin{align*}
	T_f(r,r_0)=\int\limits_{r_0}^r \dfrac{dt}{t^{2m-1}}\int_{B(t)} f^*\Omega\wedge v^{m-1}, \ (0<r_0<r<R).
\end{align*}

By  Jensen's formula, we will have
$$T_f(r,r_0)=\int\limits_{S(r)}\log\|f\|\sigma_m-\int\limits_{S(r_0)}\log\|f\|\sigma_m+O(1), \ (\text{as} \ r\to R).$$

\noindent
{\bf (e)}\ Let $H$ be a hyperplane in $\mathbb{P}^n(\mathbb{C})$ given by $H=\{a_0w_0+...+a_n+w_n=0\}$ where $(a_0,...,a_n)\neq (0,...,0)$. We set $(f,H)=\sum_{i=0}^na_if_i$. We define the proximity function by
$$m_f(r,r_0,H)=\int_{S(r)}\log\frac{\parallel f\parallel . \parallel H\parallel}{|(f,H)|}\sigma_m -\int_{S(r_0)}\log \frac{\parallel f\parallel . \parallel H\parallel}{|(f,H)|}\sigma_m$$
where $\parallel H\parallel = (\sum_{i=0}^n|a_i|^2)^{\frac{1}{2}}.$

Assume that $f(\mathbb{C}^m)\not\subset H$, then $T_f(r,r_0)=N(r,r_0,\nu_{(f,H)})+m_f(r,r_0,H)+O(1).$

As usual, by the notation ``$|| \ P$''  we mean the assertion $P$ holds for all $r \in [0,\infty)$ excluding a Borel subset $E$ of the interval $[0,\infty)$ with $\int_E\frac{1}{R-r} dr<\infty$.

\section{Proofs of Main Theorems}

Repeating the argument in [\cite{Fu1}, Proposition 4.5], we have the following

\noindent
\begin{proposition} [see \cite{Fu1}, Proposition 4.5]
	Let $F_1,\ldots,F_{n+1}$ be meromorphic functions on $B(R_0)\subset\mathbb{C}^m$ such that they are linear independent over $\mathbb{C}$. Then there exsits an admissible set $\{\alpha_i=(\alpha_{i1},\ldots,\alpha_{im})\}_{i=1}^{n+1}$ with $\alpha_{ij}\ge 0$ being integers, $|\alpha_i|=\sum_{j=1}^m|\alpha_{ij}|\le i$ for $1\le i\le n+1$ such that the generalized Wronskians $W_{\alpha_1\cdots\alpha_{n+1}}(F_1,\ldots,F_{n+1})\not\equiv 0$ where $W_{\alpha_1\cdots\alpha_{n+1}}(F_1,\ldots,F_{n+1})=\det \bigl(\mathcal{D}^{\alpha_i}F_j\bigl)_{1\le i, j \le n+1}.$
\end{proposition}

Let $L_1,\ldots,L_{n+1}$ be linear forms of $n+1$ variables and assume that they are linear independent. Let $F=(F_1:\cdots:F_{n+1}): B(R_0)\subset\mathbb{C}^m\to\mathbb{P}^n(\mathbb{C})$ be a meromorphic mapping and $(\alpha_1,\ldots,\alpha_{n+1})$ be an admissible set of $F$. Then we have following proposition.

\noindent
\begin{proposition} [see \cite{R-S1}, Proposition 3.3]
	In the above situation, set $l_0=|\alpha_1|+\cdots+|\alpha_{n+1}|$ and take $t,p$ with $0<tl_0<p<1.$ Then, for $0<r_0<R_0$ there exists a positive constant $K$ such that for $r_0<r<R<R_0,$ $$\int\limits_{S(r)}\left |z^{\alpha_1+\cdots+\alpha_{n+1}}\frac{W_{\alpha_1\cdots\alpha_{n+1}}(F_1,\ldots,F_{n+1})}{L_1(F)\cdots L_{n+1}(F)}\right|^t \sigma_m\le K\bigl(\frac{R^{2m-1}}{R-r}T_F(R,r_0)\bigl)^{p}.$$
\end{proposition}

\noindent
\begin{definition}
Let $N\ge n$ and $q\ge N+1$. Let $Q_1,\ldots,Q_q$ be hypersurfaces in $\mathbb{P}^n(\mathbb{C})$. The hypersurfaces $Q_1,\ldots,Q_q$ are said to be in $N$-subgeneral position in $\mathbb{P}^n(\mathbb{C})$ if $Q_{j_1}\cap\cdots\cap Q_{j_{N+1}}=\emptyset$ for every $1\le j_0<\cdots<j_N\le q.$
\end{definition}

If $\{Q_i\}_{i=1}^q$ is in $n$-subgeneral position then we say that it is in general position.

Let $\{Q_i\}_{i=1}^q$ be $q$ hypersurfaces in $\mathbb{P}^n(\mathbb{C})$ of the common degree $d$. Assume that each $Q_i$ is defined by a homogeneous polynomial $Q_i^*\in \mathbb{C}[x_1,\ldots,x_{n+1}].$ We regard $\mathbb{C}[x_1,\ldots,x_{n+1}]$ as a complex vector space and define $$\text{rank} \{Q_i\}_{i\in R}=\text{rank} \{Q_i^*\}_{i\in R}$$ for every subset $R\subset \{1,\ldots,q\}.$ It is easy to see that $$\text{rank} \{Q_i\}_{i\in R}=\text{rank} \{Q_i^*\}_{i\in R}\ge n+1 - \dim \bigl(\bigcap\limits_{i\in R}Q_i\bigl).$$ Hence, if $\{Q_i\}_{i=1}^q$ is in $N$-subgeneral position, by the above equality, we have $$\text{rank} \{Q_i\}_{i\in R}\ge n+1$$ for any subset $R\subset\{1,\ldots,q\}$ with $\sharp R=N+1.$

\noindent
\begin{proposition} [see \cite{Q-A}, Lemma 3.3]\label{prop3.4}
	 Let $V$ be a complex projective subvariety of $\pnc$ of dimension $k$ ($k\le n$). Let $\{Q_i\}_{i=1}^q \ (q>2N-k+1)$ be $q$ hypersurfaces of the common degree $d$ in $\mathbb{P}^n(\mathbb{C})$ located in $N$-subgeneral position with respect to $V$. Then there are positive rational constants $\omega_i \ (1\le i\le q)$ satisfying the following:
	
	(i)\ $0\le \omega_i\le 1, \ \forall i\in Q=\{1,\ldots,q\}.$
	
	(ii)\ Setting $\tilde{\omega}=\max_{i\in Q}\omega_i,$ one gets $$\sum\limits_{i=1}\limits^q \omega_i=\tilde{\omega}(q-2N+k-1)+k+1.$$
	
	(iii)\ $\frac{k+1}{2N-k+1}\le \tilde{\omega}\le \frac{k}{N}.$
	
	(iv)\ For $R\subset\{1,\ldots,q\}$ with $\sharp R=N+1,$ then $\sum_{i\in R}\omega_i\le k+1.$
	
	(v)\ Let $E_i\ge 1 \ (1\le i\le q)$ be arbitrarily given numbers. For $R\subset\{1,\ldots,q\}$ with $\sharp R$
	
	\hspace{0.5cm} $=N+1,$ there is a subset $R^o\subset R$ such that $\sharp R^o=\text{rank}\{Q_i\}_{i\in R^o}=k+1$ and $$\prod\limits_{i\in R} E_i^{\omega_i}\le \prod\limits_{i\in R^o} E_i.$$	
\end{proposition}

Let $\{Q_i\}_{i\in R}$ be a set of hypersurfaces in $\mathbb{P}^n(\mathbb{C})$ of the common degree $d$. Assume that each $Q_i$ is defined by $$\sum\limits_{I\in \mathcal{I}_d} a_{iI}x^I=0,$$ where $\mathcal{I}_d=\{I=(t_1,\ldots,t_{n+1})\in \mathbb{N}^{n+1};t_1+\cdots+t_{n+1}=d\}, \ x^I=x_1^{t_1}\cdots x_{n+1}^{t_{n+1}}$ and $(x_1:\cdots:x_{n+1})$ is homogeneous coordinates of $\mathbb{P}^n(\mathbb{C}).$

Let $f:\mathbb{C}^m\longrightarrow V\subset \mathbb{P}^n(\mathbb{C})$ be an algebraically nondegenerate meromorphic mapping with a reduced representation $f=(f_1:\cdots:f_{n+1}).$ We define $$Q_i(f)=\sum\limits_{I\in\mathcal{I}_d}a_{iI}f^I,$$ where $f^I=f_1^{t_1}\cdots f_{n+1}^{t_{n+1}}$ for $I=(t_1,\ldots,t_{n+1}).$ Then $f^* Q_i=\nu_{Q_i(f)}$ as divisors.

\noindent
\begin{lem} [see \cite{T-Q}, Lemma 3.2]
	Let $\{Q_i\}_{i\in R}$ be a set of hypersurfaces in $\mathbb{P}^n(\mathbb{C})$ of the common degree $d$ and let $f:\mathbb{C}^m\longrightarrow \mathbb{P}^n(\mathbb{C})$ be a meromorphic mapping. Assume that $\bigcap\limits_{i=1}^q Q_i\cap V=\emptyset.$ Then there exist positive constants $\alpha$ and $\beta$ such that $$\alpha\Vert f\Vert ^d\le \max\limits_{i\in R} |Q_i(f)|\le \beta\Vert f\Vert ^d.$$
\end{lem}
                                         
\noindent
\begin{lem} [see \cite{Q-A}, Lemma 4.2] 
		Let $\{Q_i\}_{i=1}^q$ be a set of hypersurfaces in $\mathbb{P}^n(\mathbb{C})$ of the common degree $d$. Then there exist $(H_V(d)-k-1)$ hypersurfaces $\{T_i\}_{i=1}^{H_V(d)-k-1}$ in $\mathbb{P}^n(\mathbb{C})$ such that for any subset $R\subset \{1,\ldots,q\}$ with $\sharp R=rank \{Q_i\}_{i\in R}=k+1$, we get $rank\{\{Q_i\}_{i\in R}\cup \{T_i\}_{i=1}^{H_V(d)-k-1}\}=H_V(d).$
\end{lem}

\noindent
\begin{lem} [see \cite{T-Q}, Lemma 3.4]
		Let $\{L_i\}_{i=1}^{H_V(d)}$ be a family of hypersurfaces in $\mathbb{P}^n(\mathbb{C})$ of the common degree $d$ and let $f$ be a meromorphic mapping of $B(R_0)\subset\mathbb{C}^m$ into $\mathbb{P}^n(\mathbb{C})$. Assume that $\{L_i\}_{i=1}^{H_V(d)}$ is linear independent. Then for every $0<r_0<r<R_0$, we have
		$$T_F(r,r_0)=dT_f(r,r_0)+O(1),$$
		where $F$ is the meromorphic mapping of $B(R_0)$ into $\mathbb{P}^n(\mathbb{C})$ with the representation $F=(L_1(f):\cdots:L_{H_V(d)}(f)).$
\end{lem}

\noindent
\textbf{Proof of Theorem \ref{1.1}}
Without loss of generality, we may assume $R_0=1.$

Assume that $f\not\equiv g$, we may choose distinct indies $i_0$ and $j_0$ such that 
$$P:=f_{i_0}g_{j_0}-f_{j_0}g_{i_0}\not\equiv 0.$$
Then $|P|\le 2\Vert f\Vert \Vert g\Vert.$

Let $Q_i \ (1\le i\le q)$ be the homogeneous polynomial in $\mathbb{C}[x_0,\ldots,x_n]$ of degree $d_i$ defining hypersurface $Q_i.$ Replacing $Q_i$ by $Q_i^{d/d_i} \ (i=1,\ldots,q)$ if necessary, we may assume that $Q_1,\ldots,Q_q$ have the same degree of $d.$

Take a $\mathbb{C}$-basis $\{[A_i] \}_{i=1}^{H_V(d)}$ of $I_d(V)$, we set 
$$F=[A_1(f):\cdots:A_{H_V(d)}(f)] \ \text{and} \ G=[A_1(g):\cdots:A_{H_V(d)}(g)].$$
Since $f$ and $g$ are nondegenerate over $I_d(V),$ $\{A_i(f)\}_{i=1}^{H_V(d)}$ and $\{A_i(g)\}_{i=1}^{H_V(d)}$ are linear independent over $\mathbb{C}.$ Then there exist an admissible sets $\{\alpha_1,\ldots,\alpha_{H_V(d)} \}$ and $\{\beta_1,\ldots,\beta_{H_V(d)} \}$ such that 
$$W_{\alpha_1\cdots\alpha_{H_V(d)}}(F)=\det \bigl(\mathcal{D}^{\alpha_i}A_j(f)\bigl)_{1\le i, j \le H_V(d)}\not\equiv 0$$
and 
$$W_{\beta_1\cdots\beta_{H_V(d)}}(G)=\det \bigl(\mathcal{D}^{\alpha_i}A_j(g)\bigl)_{1\le i, j \le H_V(d)}\not\equiv 0.$$

Let $z$ be a fixed point. Then there exists $R\subset Q=\{1,\ldots,q\}$ with $\sharp R=N+1$ such that $|Q_i(f)(z)|\le|Q_j(f)(z)|, \ \forall i\in R, \ j\notin R$. Since $\bigcap\limits_{i\in R} Q_i\bigcap V=\emptyset$, there exist positive constants $\alpha$ and $\beta$ such that 
$$\alpha\Vert f\Vert ^d\le \max\limits_{i\in R} |Q_i(f)|\le \beta\Vert f\Vert ^d.$$

On the other hand, there exists a subset $R^0\subset R$ such that $\sharp R^0=k+1$ and $R^0$ satisfies Proposition \ref{prop3.4} v) with respect to number $\left\{\frac{\beta \Vert f\Vert^d}{Q_i(f)(z)}\right\}_{i=1}^q.$ We get 
$$\mathrm{rank}\{\{Q_i\}_{i\in R^0} \cup \{T_i\}_{i=1}^{H_V(d)-k-1} \}=H_V(d),$$
where $\{\{T_i\}_{i=1}^{H_V(d)-k-1} \}$ are hypersurfaces in $\mathbb{P}(\mathbb{C})$ with respect to $V.$
For $R^0=\{r^0_1,\ldots,r^0_{k+1} \}$, we set
$$W_{R^0}=\det \bigl(D^{\alpha_i} Q_{r^0_v}(f) \ (1\le v\le k+1), D^{\alpha_i} T_l(f) \ (1\le l\le H_V(d)-k-1)\bigl)_{1\le i\le H_V(d)}.$$
Then there exists a nonzero constant $C_{R^0}$ such that $W_{R^0}=C_{R^0}W_{\alpha_1\cdots\alpha_{H_V(d)}}(F).$

For positive rational constants $\omega_i \ (1\le i\le q),$ we have 
\begin{align*}
\frac{\Vert f(z)\Vert ^{d(\sum\limits_{i=1}^q \omega_i)}|W_{\alpha_1\cdots\alpha_{H_V(d)}}(F)(z)|}{|Q_1^{\omega_1}(f)(z)\cdots Q_q^{\omega_q} (f)(z)|}
&\le \frac{|W_{\alpha_1\cdots\alpha_{H_V(d)}}(F)(z)|}{\alpha^{q-N-1}\beta^{N+1}}\prod\limits_{i\in R}\bigl( \frac{\beta \Vert f(z)\Vert^d}{|Q_i(f)(z)|}\bigl)^{\omega_i}\\  
&\le K \frac{|W_{\alpha_1\cdots\alpha_{H_V(d)}}(F)(z)|\Vert f(z)\Vert ^{d(k+1)}}{\prod\limits_{i\in R^0} |Q_i(f)(z)|}\\
&\le K_0  \frac{|W_{R^0}(z)|\Vert f(z)\Vert ^{dH_V(d)}}{\prod\limits_{i\in R^0} |Q_i(f)(z)|\prod\limits_{i=1}^{H_V(d)-k-1} |T_i(f)(z)|}.
\end{align*}
where $K$ and $K_0$ are positive constants.

Hence 
$$\frac{\Vert f(z)\Vert ^{d(\sum\limits_{i=1}^q \omega_i-H_V(d))}|W_{\alpha_1\cdots\alpha_{H_V(d)}}(F)(z)|}{|Q_1^{\omega_1}(f)(z)\cdots Q_q^{\omega_q} (f)(z)|}
\le K_0  \frac{|W_{R^0}(z)|}{\prod\limits_{i\in R^0} |Q_i(f)(z)|\prod_{i=1}^{H_V(d)-k-1} |T_i(f)(z)|}.$$

Put $S_{R^0}=\frac{W_{R^0}}{\prod\limits_{i\in R^0} Q_i(f)\prod_{i=1}^{H_V(d)-k-1} T_i(f)}$. Then for each $z\in\mathbb{C}^m$, we get
$$\frac{\Vert f(z)\Vert ^{d(\sum\limits_{i=1}^q \omega_i-H_V(d))}|W_{\alpha_1\cdots\alpha_{H_V(d)}}(F)(z)|}{|Q_1^{\omega_1}(f)(z)\cdots Q_q^{\omega_q} (f)(z)|}
\le K_0|S_{R^0}(z)|.$$

Letting $\phi=z^{\alpha_1+\cdots+\alpha_{H_V(d)}}\frac{W_{\alpha_1\cdots\alpha_{H_V(d)}}(F)}{Q_1^{\omega_1}(f)\cdots Q_q^{\omega_q}(f)}$ and $\psi=z^{\beta_1+\cdots+\beta_{H_V(d)}}\frac{W_{\beta_1\cdots\beta_{H_V(d)}}(G)}{Q_1^{\omega_1}(g)\cdots Q_q^{\omega_q}(g)},$ we have
$$\sum\limits_{i=1}^q\omega_i \nu_{Q_i(f)}(z)-\nu_{W_{\alpha_1\cdots\alpha_{H_V(d)}}(F)}(z)\le \sum\limits_{i=1}^q\omega_i\min\{H_V(d)-1,\nu_{Q_i(f)}(z)\},$$
$$\sum\limits_{i=1}^q\omega_i \nu_{Q_i(g)}(z)-\nu_{W_{\beta_1\cdots\beta_{H_V(d)}}(G)}(z)\le \sum\limits_{i=1}^q\omega_i\min\{H_V(d)-1,\nu_{Q_i(g)}(z)\}.$$
Indeed, let $z$ be a zero of some $Q_i(f)$ and $z\notin I(f)=\{f_0=\cdots=f_n=0\}$. Since $\{Q_i\}_{i=1}^q$ is in $N$-subgeneral position, $z$ is not zero of more than $N$ functions $Q_i(f)$. Without loss of generality, we may assume that $z$ is not zero of $Q_i(f)$ for each $i>N$. Put $R=\{1,\ldots, N+1\}$ and choose $R^1\subset R$ such that $\sharp R^1=\text{rank}\{Q_i\}_{i\in R^1}=k+1$ satisfied Proposition \ref{prop3.4} v) with respect to numbers $\{e^{\max \{\nu_{Q_i(f)}(z)-H_V(d)+1,0\}}\}_{i=1}^q$. Then we have 
$$\sum\limits_{i\in R} \omega_i \max \{\nu_{Q_i(f)}(z)-H_V(d)+1,0\}\le \sum\limits_{i\in R^1}\max \{\nu_{Q_i(f)}(z)-H_V(d)+1,0\}.$$
This implies that
\begin{align*}
\nu_{W_{\alpha_1\cdots\alpha_{H_V(d)}}(F)}(z)& \ge  \sum\limits_{i\in R^1}\max \{\nu_{Q_i(f)}(z)-H_V(d)+1,0\}\\
&\ge \sum\limits_{i\in R} \omega_i \max \{\nu_{Q_i(f)}(z)-H_V(d)+1,0\}.
\end{align*}
Hence 
\begin{align*}
\sum\limits_{i=1}^q\omega_i\nu_{Q_i(f)}(z)-&\nu_{W_{\alpha_1\cdots\alpha_{H_V(d)}}(F)}(z)=\sum\limits_{i\in R}\omega_i \nu_{Q_i(f)}(z)-\nu_{W_{\alpha_1\cdots\alpha_{H_V(d)}}(F)}(z)\\
&=\sum\limits_{i\in R}\omega_i \min\{\nu_{Q_i(f)}(z), H_V(d)-1\}\\
&\quad +\sum\limits_{i\in R}\omega_i \max\{\nu_{Q_i(f)}(z) - H_V(d)+1,0\}-\nu_{W_{\alpha_1\cdots\alpha_{H_V(d)}}(F)}(z)\\
&\le \sum\limits_{i\in R}\omega_i \min\{\nu_{Q_i(f)}(z), H_V(d)-1\}.
\end{align*}
This proves the above inequalities.

On the other hand, since $P=0$ on $\bigcup_{i=1}^q\bigl(Zero Q_i(f)\cup Zero Q_i(g)\bigl)$ and $\{Q_i\}_{i=1}^q$ is in $N$-subgeneral position, there are at most $N$ functions $Q_i(f)$ vanishing at each point of $\bigcup_{i=1}^q\bigl(Zero Q_i(f)\cup Zero Q_i(g)\bigl)$. Therefore, we have,
$$\nu_P(z)\ge \frac{1}{N}\sum\limits_{i=1}^q\min\{1,\nu_{Q_i(f)}(z)\}$$
It implies that 
$$\nu_P(z)\ge \frac{1}{N(H_V(d)-1)}\sum\limits_{i=1}^q\min\{H_V(d)-1,\nu_{Q_i(f)}(z)\}.$$
Similarly, we also get 
$$\nu_P(z)\ge \frac{1}{N(H_V(d)-1)}\sum\limits_{i=1}^q\min\{H_V(d)-1,\nu_{Q_i(g)}(z)\}.$$

Therefore, we can see that
$$\sum\limits_{i=1}^q\omega_i \nu_{Q_i(f)}(z)-\nu_{W_{\alpha_1\cdots\alpha_{H_V(d)}}(F)}(z)\le \tilde{\omega}N(H_V(d)-1)\nu_P(z),$$
$$\sum\limits_{i=1}^q\omega_i \nu_{Q_i(g)}(z)-\nu_{W_{\beta_1\cdots\beta_{H_V(d)}}(G)}(z)\le\tilde{\omega}N(H_V(d)-1)\nu_P(z).$$
Hence, $\phi P^{\tilde{\omega}N(H_V(d)-1)}$ and $\psi P^{\tilde{\omega}N(H_V(d)-1)}$ are both holomorphic on $B(1).$
If we let 
$$t=\frac{\rho}{d\tilde{\omega}\bigl(q-2N+k-1-(H_V(d)-k-1)/\tilde{\omega}-2N(H_V(d)-1)/d\bigl)},$$ 
and 
$$v=\log|\phi\psi P^{2\tilde{\omega}N(H_V(d)-1)}|,$$
then $u=tv$ is plurisubharmonic on $B(1)$, and hence
$$2t\bigl(\sum\limits_{s=1}^{H_V(d)}|\alpha_s| \bigl)< \frac{2\rho}{\rho H_V(d)(H_V(d)-1)}.\frac{H_V(d)(H_V(d)-1)}{2}=1,$$
for
\begin{align*}
	q&>\frac{(2N-k+1)H_V(d)}{k+1}\bigl(1+\frac{\rho(H_V(d)-1)}{d}\bigl)+\frac{2N(H_V(d)-1)}{d}\\
	&\ge 2N-k+1+\frac{(H_V(d)-k-1)}{\tilde{\omega}}+\frac{2N(H_V(d)-1)}{d}+\frac{\rho(2N-k+1)H_V(d)(H_V(d)-1)}{(k+1)d}.
\end{align*}
Then we may choose a positive number $p$ such that 
$$0<2t\bigl(\sum\limits_{s=1}^{H_V(d)}|\alpha_s| \bigl)<p<1.$$

We write the given K\"{a}hler metric form as 
$$\omega=\frac{\sqrt{-1}}{2\pi}\sum\limits_{i,j}h_{i\bar{j}}dz_i\wedge d\bar{z}_j.$$
From the assumption that both $f$ and $g$ satisfy condition $(C_\rho)$, there are continuous plurisubharmonic functions $u_1, u_2$ on $B(1)$ such that
$$e^{u_1}\text{det}(h_{i\bar{j}})^{\frac{1}{2}}\le \Vert f\Vert ^\rho,$$ $$e^{u_2}\text{det}(h_{i\bar{j}})^{\frac{1}{2}}\le \Vert g\Vert ^\rho.$$ 
Therefore 
\begin{align*}
e^{u+u_1+u_2}&\text{det}(h_{i\bar{j}})=e^{tv+u_1+u_2}\text{det}(h_{i\bar{j}})\\
&\le e^{tv}\Vert f\Vert ^\rho \Vert g\Vert ^\rho=|\phi|^t|\psi|^t|P|^{2t\tilde{\omega}N(H_V(d)-1)}\Vert f\Vert ^\rho\Vert g\Vert ^\rho\\
&\le C |\phi|^t|\psi|^t\Vert f\Vert ^{\rho+2t\tilde{\omega}N(H_V(d)-1)}\Vert g\Vert ^{\rho+2t\tilde{\omega}N(H_V(d)-1)}\\
&\le C |\phi|^t|\psi|^t\Vert f\Vert ^{td\tilde{\omega}(q-2N+k-1-(H_V(d)-k-1)/\tilde{\omega})}\Vert g\Vert ^{td\tilde{\omega}(q-2N+k-1-(H_V(d)-k-1)/\tilde{\omega})}\\
& \le C |\phi|^t|\psi|^t\Vert f\Vert ^{td(\sum_{i=1}^q\omega_i-H_V(d))}\Vert g\Vert ^{td(\sum_{i=1}^q\omega_i-H_V(d))}.
\end{align*}
for $C$ is a positive constant. Note that the volume form on $B(1)$ is given by 
$$dV:=c_m\text{det}(h_{i\bar{j}})v_m;$$
therefore, 
$$\int\limits_{B(1)} e^{u+u_1+u_2}dV\le C\int\limits_{B(1)} |\phi|^t|\psi|^t\Vert f\Vert ^{td(\sum_{i=1}^q\omega_i-H_V(d))}\Vert g\Vert ^{td(\sum_{i=1}^q\omega_i-H_V(d))} v_m.$$
Thus, by the H\"{o}lder inequality and by noticing that 
$$v_m=(dd^c\Vert z\Vert^2)^m=2m\Vert z\Vert^{2m-1}\sigma_m\wedge d\Vert z\Vert,$$
we obtain 
\begin{align*}
\int\limits_{B(1)} e^{u+u_1+u_2}dV&\le C\bigl (\int\limits_{B(1)} |\phi|^{2t}\Vert f\Vert ^{2td(\sum_{i=1}^q\omega_i-H_V(d))} v_m \bigl)^{\frac{1}{2}}\\
&\hspace{3cm}\times\bigl(\int\limits_{B(1)} |\psi|^{2t}\Vert g\Vert ^{2td(\sum_{i=1}^q\omega_i-H_V(d))} v_m \bigl)^{\frac{1}{2}}\\
&\le C\bigl(2m\int\limits_0\limits^1 r^{2m-1}\bigl(\int\limits_{S(r)} |\phi|^{2t}\Vert f\Vert ^{2td(\sum_{i=1}^q\omega_i-H_V(d))} \sigma_m\bigl)dr\bigl)^{\frac{1}{2}}\\
&\hspace{3cm}\times \bigl(2m\int\limits_0\limits^1 r^{2m-1}\bigl(\int\limits_{S(r)} |\psi|^{2t}\Vert g\Vert ^{2td(\sum_{i=1}^q\omega_i-H_V(d))}\sigma_m\bigl)dr\bigl)^{\frac{1}{2}}.
\end{align*}

We note that
\begin{align*}
\int\limits_{S(r)} |\phi|^{2t}&\Vert f\Vert ^{2td(\sum_{i=1}^q\omega_i-H_V(d))} \sigma_m\\
&=\int\limits_{S(r)}\left|z^{\alpha_1+\cdots+\alpha_{H_V(d)}}\frac{W_{\alpha_1\cdots\alpha_{H_V(d)}}(F)}{Q_1^{\omega_1}(f)\cdots Q_q^{\omega_q}(f)} \right|^{2t}\bigl(\Vert f\Vert ^{d(\sum_{i=1}^q\omega_i-H_V(d))}\bigl)^{2t} \sigma_m \\
&\le \int\limits_{S(r)}\left|z^{\alpha_1+\cdots+\alpha_{H_V(d)}}K_0 S_{R^0}\right|^{2t}\sigma_m.
\end{align*}

On the other hand, for $p$ satisfying $0<2t(\sum\limits_{s=1}^{H_V(d)}|\alpha_s|)<p<1$ and for $0<r_0<r<R<1$, we get
$$\int\limits_{S(r)}\left|z^{\alpha_1+\cdots+\alpha_{H_V(d)}}K_0 S_{R^0}\right|^{2t}\sigma_m\le K_1\bigl(\frac{R^{2m-1}}{R-r} dT_f(R,r_0) \bigl)^p.$$
outside a subset $E\subset [0,1]$ such that $\int\limits_{E}\frac{1}{1-r}dr\le +\infty.$ Choosing $R=r+\frac{1-r}{eT_f(r,r_0)},$ we have
$$T_f(R,r_0)\le 2T_f(r,r_0),$$

Hence, the above inequality implies that 
$$\int\limits_{S(r)}\left|z^{\alpha_1+\cdots+\alpha_{H_V(d)}}K_0 S_{R^0}\right|^{2t}\sigma_m\le\frac{K'}{(1-r)^p}(T_f(r,r_0))^{2p}
\le \frac{K'}{(1-r)^p}(\log\frac{1}{1-r})^{2p},$$
since $\lim  \limits_{r\to R_0}\sup\frac{T_f(r,r_0)+T_g(r,r_0)}{\log (1/(R_0-r))}<\infty.$

It implies that
$$2m\int\limits_0\limits^1 r^{2m-1}\bigl(\int\limits_{S(r)} |\phi|^{2t}\Vert f\Vert ^{2td(\sum_{i=1}^q-H_V(d))} \sigma_m\bigl)dr\le 2m\int\limits_0\limits^1 r^{2m-1}\frac{K'}{(1-r)^p}(\log\frac{1}{1-r})^{2p} dr <\infty.$$
Similarly,
$$2m\int\limits_0\limits^1 r^{2m-1}\bigl(\int\limits_{S(r)} |\psi|^{2t}\Vert g\Vert ^{2td(\sum_{i=1}^q-H_V(d))} \sigma_m\bigl)dr\le 2m\int\limits_0\limits^1 r^{2m-1}\frac{K'}{(1-r)^p}(\log\frac{1}{1-r})^{2p} dr <\infty.$$

Hence, we conclude that
$$\int\limits_{B(1)} e^{u+u_1+u_2}dV<\infty,$$
which contradicts the result of Yau \cite{Y} and Karp\cite{K}. 

This completes the proof.\hfill$\square$

\vskip0.2cm
In order to prove Theorem \ref{1.3}, we need the following.
\begin{lem} [see \cite{No2}, Lemma 3.3]
	Let $\{H_i\}_{i=1}^q \ (q>2N-n+1)$ be $q$ hyperplanes in $\mathbb{P}^n(\mathbb{C})$ located in $N$-subgeneral position. Then there are positive rational constants $\omega_i \ (1\le i\le q)$ satisfying the following:
	
	(i)\ $0\le \omega_i\le 1, \ \forall i\in Q=\{1,\ldots,q\}.$
	
	(ii)\ Setting $\tilde{\omega}=\max_{i\in Q}\omega_i,$ one gets $$\sum\limits_{i=1}\limits^q \omega_i=\tilde{\omega}(q-2N+n-1)+n+1.$$
	
	(iii)\ $\frac{n+1}{2N-n+1}\le \tilde{\omega}\le \frac{n}{N}.$
	
	(iv)\ For $R\in\{1,\ldots,q\}$ with $\sharp R=N+1,$ then $\sum_{i\in R}\omega_i\le n+1.$
	
	(v)\ Let $E_i\ge 1 \ (1\le i\le q)$ be arbitrarily given numbers. For $R\in\{1,\ldots,q\}$ with $\sharp R$
	
	\hspace{0.5cm} $=N+1,$ there is a subset $R^o\subset R$ such that $\sharp R^o=\text{rank}\{H_i\}_{i\in R^o}=n+1$ and $$\prod\limits_{i\in R} E_i^{\omega_i}\le \prod\limits_{i\in R^o} E_i.$$
\end{lem}

\noindent
\textbf{Proof of Theorem\ref{1.3}.}\
By using the universal covering if necessary, it suffices for us to prove the theorem in the case where $M$ is the ball $B(1)$ of $\mathbb{C}^m.$

We assume that each $H_i$ is given by $$H_i=\{(z_0:\cdots:z_n)\in\mathbb{P}^n(\mathbb{C}); a_{i0}z_0+\cdots+a_{in}z_n=0\},$$
where $a_{ij}\in\mathbb{C}$ and $\sum\limits_{j=0}^n|a_{ij}|^2=1.$ We set $$H_i(f)=a_{i0}f_0+\cdots+a_{in}f_n.$$
Hence, it is easy to see that $|H_i(f)|\le \lVert f\lVert$ for each $i.$

Since $f$ is differential nondegenerate, there exists a set of indices $\{j_1,\ldots,j_n\}\subset \{1,\ldots,m\}$ such that

\[
W(f)=W(f_0,\ldots,f_n):=\left |\begin{array}{cccc}
f_0 &f_1 &\cdots &f_n\\
\frac{\partial f_0}{\partial z_{j1}} &\frac{\partial f_1}{\partial z_{j1}} &\cdots &\frac{\partial f_n}{\partial z_{j1}}\\
\vdots &\vdots &\ddots &\vdots\\
\frac{\partial f_0}{\partial z_{jn}} &\frac{\partial f_1}{\partial z_{jn}} &\cdots &\frac{\partial f_n}{\partial z_{jn}}\\
\end{array}
\right |\not\equiv 0
\]
Without loss of generality, we may assume that the set of indices $\{j_1,\ldots,j_n\}$ is $\{1,\ldots,n\}.$

For each $R^0=\{r^0_1,\ldots,r^0_{n+1}\}\subset \{1,\ldots,q\}$ with $\mathrm{rank}\{H_i\}_{i\in R^0}=\sharp R^0=n+1,$ we set
$$
W_{R^0}(f)=\left |\begin{array}{cccc}
H_{r^0_1}(f) &H_{r^0_2}(f) &\cdots &H_{r^0_{n+1}}(f)\\
\frac{\partial H_{r^0_1}(f)}{\partial z_{1}} &\frac{\partial H_{r^0_2}(f)}{\partial z_{1}} &\cdots &\frac{\partial H_{r^0_{n+1}}(f)}{\partial z_{1}}\\
\vdots &\vdots &\ddots &\vdots\\
\frac{\partial H_{r^0_1}(f)}{\partial z_{n}} &\frac{\partial H_{r^0_2}(f)}{\partial z_{n}} &\cdots &\frac{\partial H_{r^0_{n+1}}(f)}{\partial z_{n}}\\
\end{array}
\right |.
$$
Since $\text{rank}\{H_{r^0_v}(1\le v\le n+1)\}=n+1,$ there exist a nonzero constant $C_{R^0}$ such that $W_{R^0}(f)=C_{R^0}.W(f).$

We denote by $\mathcal{R}^0$ the family of all subsets $R^0$ of $\{1,\ldots,q\}$ with $\text{rank}\{H_i\}_{i\in R^0}=\sharp R^0=n+1.$

Let $z$ be a fixed point. Take a subset $R\subset Q$ with $\sharp R=N+1,$ such that 
	$$|H_i(f)(z)|\le |H_j(f)(z)|, \forall i\in R, j\notin R.$$
Since $\cap_{i\in R}H_i=\emptyset,$ there exist positive constants $\alpha$ and $\beta$ such that 
$$\alpha\Vert f\Vert\le\max\limits_{i\in R}|H_i(f)|\le \beta\Vert f\Vert.$$

On the other hand, we may choose $R^0\subset R$ such that $R^0\in\mathcal{R}^0$ and $R^0$ satisfies Lemma 3.8 iv) with respect to numbers  $\{\frac{\beta \Vert f(z)\Vert}{|H_i(f)(z)|}\}_{i=1}^q.$ For positive rational constant $\omega_i \ (1\le i\le q),$ we get 
\begin{align*}
\frac{\Vert f(z)\Vert ^{(\sum\limits_{i=1}^q \omega_i)}|W(f)(z)|}{|H_1^{\omega_1}(f)(z)\cdots H_q^{\omega_q} (f)(z)|}
&\le \frac{|W(f)(z)|}{\alpha^{q-N-1}\beta^{N+1}}\prod\limits_{i\in R}\bigl( \frac{\beta \Vert f(z)\Vert^d}{|H_i(f)(z)|}\bigl)^{\omega_i}\\  
&\le K_0 \frac{|W_{R^0}(f)(z)|\Vert f(z)\Vert ^{(n+1)}}{\prod\limits_{i\in R^0} |H_i(f)(z)|},
\end{align*}
where $K_0$ is a positive constant.

Put $S_{R^0}=\frac{W_{R^0}(f)}{\prod\limits_{i\in R^0} H_i(f)}$. Then for each $z\in\mathbb{C}^m$, we get
$$\frac{\Vert f(z)\Vert ^{(\sum\limits_{i=1}^q \omega_i-n-1)}|W(f)(z)|}{|H_1^{\omega_1}(f)(z)\cdots H_q^{\omega_q} (f)(z)|}\le K_0.|S_{R^0}(z)|.$$
Let $t=\frac{\rho}{\tilde{\omega}(q-2N+n-3)}$ and $v=\log|\phi\psi P^{2\tilde{\omega}}|,$ where $\phi=(\prod\limits_{i=1}^{n+1}z_i).\frac{W(f)}{H_1^{\omega_1}(f)\cdots H_q^{\omega_q}(f)}$ and $\psi=(\prod\limits_{i=1}^{n+1}z_i).\frac{W(g)}{H_1^{\omega_1}(g)\cdots H_q^{\omega_q}(g)}.$

Repeating the argument in proof of Theorem \ref{1.1}, we get $u=tv$ is plurisubharmonic on $B(1)$ and 
$$2t.(n+1)< \frac{2\rho}{2\rho(n+1)}.(n+1)=1,$$
for
$q>2N-n+3+2\rho(2N-n+1)\ge 2N-n+3+\frac{2\rho(n+1)}{\tilde{\omega}}.$

Then we may choose a positive number $p$ such that 
$$0<2t(n+1)<p<1.$$
From the assumption that both $f$ and $g$ satisfy the condition $(C_\rho)$, there are continuous plurisubharmonic functions $u_1, u_2$ on $B(1)$ such that
$$e^{u_1}\text{det}(h_{i\bar{j}})^{\frac{1}{2}}\le \Vert f\Vert ^\rho,$$ $$e^{u_2}\text{det}(h_{i\bar{j}})^{\frac{1}{2}}\le \Vert g\Vert ^\rho.$$ 
Therefore 
\begin{align*}
e^{u+u_1+u_2}\text{det}(h_{i\bar{j}})&=e^{tv+u_1+u_2}\text{det}(h_{i\bar{j}})\\
&\le e^{tv}\Vert f\Vert ^\rho \Vert g\Vert ^\rho=|\phi|^t|\psi|^t|P|^{2t\tilde{\omega}}\Vert f\Vert ^\rho\Vert g\Vert ^\rho\\
&\le C |\phi|^t|\psi|^t\Vert f\Vert ^{\rho+2t\tilde{\omega}}\Vert g\Vert ^{\rho+2t\tilde{\omega}}\\
&\le C |\phi|^t|\psi|^t\Vert f\Vert ^{t\tilde{\omega}(q-2N+n-1}\Vert g\Vert ^{t\tilde{\omega}(q-2N+n-1}\\
& \le C |\phi|^t|\psi|^t\Vert f\Vert ^{t(\sum_{i=1}^q\omega_i-n-1)}\Vert g\Vert ^{t(\sum_{i=1}^q\omega_i-n-1)},
\end{align*}
where $C$ is a positive constant. Note that the volume form on $B(1)$ is given by 
$$dV:=c_m\text{det}(h_{i\bar{j}})v_m;$$
therefore, 
$$\int\limits_{B(1)} e^{u+u_1+u_2}dV\le C\int\limits_{B(1)} |\phi|^t|\psi|^t\Vert f\Vert ^{t(\sum_{i=1}^q\omega_i-n-1)}\Vert g\Vert ^{t(\sum_{i=1}^q\omega_i-n-1)} v_m.$$
Thus, by the H\"{o}lder inequality and by noticing that $$v_m=(dd^c\Vert z\Vert^2)^m=2m\Vert z\Vert^{2m-1}\sigma_m\wedge d\Vert z\Vert,$$
we obtain 
\begin{align*}
\int\limits_{B(1)} e^{u+u_1+u_2}dV&\le C\bigl(\int\limits_{B(1)} |\phi|^{2t}\Vert f\Vert ^{2t(\sum_{i=1}^q\omega_i-n-1)} v_m \bigl)^{\frac{1}{2}}.\bigl(\int\limits_{B(1)} |\psi|^{2t}\Vert g\Vert ^{2t(\sum_{i=1}^q\omega_i-n-1)} v_m \bigl)^{\frac{1}{2}}\\
&\le C\bigl(2m\int\limits_0\limits^1 r^{2m-1}\bigl(\int\limits_{S(r)} |\phi|^{2t}\Vert f\Vert ^{2t(\sum_{i=1}^q\omega_i-n-1)} \sigma_m\bigl)dr\bigl)^{\frac{1}{2}}\\
&\hspace{3cm}\times \bigl(2m\int\limits_0\limits^1 r^{2m-1}\bigl(\int\limits_{S(r)} |\psi|^{2t}\Vert g\Vert ^{2t(\sum_{i=1}^q\omega_i-n-1)}\sigma_m\bigl)dr\bigl)^{\frac{1}{2}}.
\end{align*}

We may see that
\begin{align*}
\int\limits_{S(r)} |\phi|^{2t}\Vert f\Vert ^{2t(\sum_{i=1}^q\omega_i-n-1)} \sigma_m&=\int\limits_{S(r)}\left|(\prod\limits_{i=1}^{n+1}z_i).\frac{W(f)}{H_1^{\omega_1}(f)\cdots H_q^{\omega_q}(f)} \right|^{2t}\bigl(\Vert f\Vert ^{(\sum_{i=1}^q\omega_i-n-1)}\bigl)^{2t} \sigma_m \\
&\le \int\limits_{S(r)}\left|(\prod\limits_{i=1}^{n+1}z_i)K_0 S_{R^0}\right|^{2t}\sigma_m.
\end{align*}
On the other hand, for $p$ satisfy $0<2t(n+1)<p<1$ and for $0<r_0<r<R<1$, we get $$\int\limits_{S(r)}\left|(\prod\limits_{i=1}^{n+1}z_i)K_0 S_{R^0}\right|^{2t}\sigma_m\le K_1\bigl(\frac{R^{2m-1}}{R-r} dT_f(R,r_0) \bigl)^p.$$
outside a set $E\subset[0,1]$ such that $\int\limits_{E}\frac{1}{1-r}dr\le +\infty.$ Choosing $R=r+\frac{1-r}{eT_f(r,r_0)},$ we have
$$T_f(R,r_0)\le 2T_f(r,r_0),$$

Hence, the above inequality implies that 
$$\int\limits_{S(r)}\left|(\prod\limits_{i=1}^{n+1}z_i)K_0 S_{R^0}\right|^{2t}\sigma_m\le\frac{K'}{(1-r)^p}(T_f(r,r_0))^{2p}
\le \frac{K'}{(1-r)^p}(\log\frac{1}{1-r})^{2p}$$
since $\lim \limits_{r\to R_0} \sup\frac{T_f(r,r_0)+T_g(r,r_0)}{\log (1/(R_0-r))}<\infty.$

It implies that
$$2m\int\limits_0\limits^1 r^{2m-1}\bigl(\int\limits_{S(r)} |\phi|^{2t}\Vert f\Vert ^{2t(\sum_{i=1}^q\omega_i-n-1)} \sigma_m\bigl)dr\le 2m\int\limits_0\limits^1 r^{2m-1}\frac{K'}{(1-r)^p}(\log\frac{1}{1-r})^{2p} dr <\infty.$$
Similarly, we have
$$2m\int\limits_0\limits^1 r^{2m-1}\bigl(\int\limits_{S(r)} |\psi|^{2t}\Vert g\Vert ^{2t(\sum_{i=1}^q\omega_i-n-1)} \sigma_m\bigl)dr\le 2m\int\limits_0\limits^1 r^{2m-1}\frac{K'}{(1-r)^p}(\log\frac{1}{1-r})^{2p} dr <\infty.$$

Hence, we conclude that
$$\int\limits_{B(1)} e^{u+u_1+u_2}dV<\infty,$$
which contradicts the result of Yau\cite{Y} and Karp\cite{K}. 

This completes the proof.
\hfill$\square$

\end{document}